\pdfoutput=0
\documentclass{article}
 \usepackage{amsmath}
 \usepackage{amsfonts}
 \usepackage{amsthm}
 \usepackage{amssymb}
 \usepackage{array}
 \usepackage{multirow}
\usepackage{graphicx}
 \usepackage{hyperref}
 \usepackage{algorithm}
 \usepackage{algpseudocode}
 \usepackage{mathtools}
 \usepackage{authblk}
\DeclarePairedDelimiter\floor{\lfloor}{\rfloor} 

\newcommand{\be}[1]{\begin{equation}\label{#1}}
\newcommand{\benon}{\begin{equation*}}  
\newcommand{\bemuln}[1]{\begin{multline}\label{#1}}
\newcommand{\bemul}{\begin{multline*}}
\newcommand{\bee}{\begin{eqnarray*}}
\newcommand{\eee}{\end{eqnarray*}}
\newcommand{\been}[1]{\begin{eqnarray}\label{#1}}
\newcommand{\eeen}{\end{eqnarray}}
\newcommand{\began}[1]{\begin{gather}\label{#1}}
\newcommand{\bega}{\begin{gather*}}
\newcommand{\bealn}[1]{\begin{align}\label{#1}}
\newcommand{\beal}{\begin{align*}}
\newcommand{\bealatn}[2]{\begin{alignat}{#1}\label{#2}}
\newcommand{\bealat}{\begin{alignat*}}
\newcommand{\bexalatn}[1]{\begin{xalignat}\label{#1}}
\newcommand{\bexalat}{\begin{xalignat*}}

\newcommand{\mbb}{\mathbb}

\newtheorem{thm}{Theorem}[section]

\def\by{{\mathbf y}}

\def\texitem#1{\par\smallskip\noindent\hangindent 25pt
               \hbox to 25pt {\hss #1 ~}\ignorespaces}

\newcommand{\scrI}{\mathcal{I}}

\newcommand{\scrS}{\mathcal{S}}

\newcommand{\balpha}{\boldsymbol{\alpha}}
\newcommand{\bbeta}{\boldsymbol{\beta}}

\newcommand{\bepsilon}{\boldsymbol{\epsilon}}

\newcommand{\btheta}{\boldsymbol{\theta}}

\usepackage{cite}

\usepackage{xcolor}

\newcommand\hl[1]{%
  \bgroup
  \hskip0pt\color{red!80!black}%
  #1%
  \egroup
}

\newcommand\hb[1]{%
	\bgroup
	\hskip0pt\color{blue!80!black}%
	#1%
	\egroup
}

\begin{document}

\title{Learning from Past Bids to Participate Strategically in Day-Ahead
  Electricity Markets$^*$ \thanks{* Research partially supported by the
    NSF under grants DMS-1664644, CNS-1645681, CCF-1527292, AitF-9500307423, and
    IIS-1237022, by the ARO under grant W911NF-12-1-0390, by the ONR
    under grant N00014-16-1-2832, and by the Sloan Foundation under
    grant G-2017-9723.}}

\date{}

\author[1]{Ruidi Chen}
\author[1,2,3]{Ioannis~Ch.~Paschalidis }
\author[1]{Michael~C.~Caramanis }
\author[1]{Panagiotis~Andrianesis }
\affil[1]{Division of Systems Engineering, Boston University}
\affil[2]{Department of
	Electrical and Computer Engineering}
\affil[3]{Department of
	Biomedical Engineering}

  
\maketitle

\begin{abstract}
  We consider the process of bidding by electricity suppliers in a
  day-ahead market context where each supplier bids a linear
  non-decreasing function of her generating capacity with the goal of
  maximizing her individual profit given other competing suppliers'
  bids. Based on the submitted bids, the market operator schedules
  suppliers to meet demand during each hour and determines hourly market
  clearing prices. Eventually, this game-theoretic process reaches a
  Nash equilibrium when no supplier is motivated to modify her
  bid. However, solving the individual profit maximization problem
  requires information of rivals' bids, which are typically not
  available. To address this issue, we develop an inverse optimization
  approach for estimating rivals' production cost functions given
  historical market clearing prices and production levels. We then use
  these functions to bid strategically and compute Nash equilibrium
  bids. We present numerical experiments illustrating our methodology,
  showing good agreement between bids based on the estimated production
  cost functions with the bids based on the true cost functions. We
  discuss an extension of our approach that takes into account network
  congestion resulting in location-dependent prices.
\end{abstract}


\section{Introduction}

In the past several decades that followed the seminal
work on spot market pricing~\cite{mcc1988spot}, the electricity industry
has evolved from vertical integrated regulated monopolies to competitive
supply and demand market participants with equal access to a regulated
transmission and distribution network.  Nevertheless, due to special
features of the power industry, including a limited number of producers
(electricity suppliers), large capital investments that introduce
barriers to entry, and congestion caused by occasionally binding
transmission constraints, the electricity market is characterized by
oligopolistic conditions~\cite{david2000strategic}.  Under perfect
competition, suppliers would bid their marginal costs, a necessary
condition for social welfare and efficiency maximization. In an
imperfect oligopolistic energy market setting, however, suppliers can
exploit market manipulation opportunities to increase their profits by
bidding above their marginal cost.

The investigation of such behavior, referred to as {\em strategic
  bidding}, is of dual interest.  First, to market participants (mainly
suppliers), who are interested in devising optimal bidding strategies
that would allow them to ``outsmart'' competitors and realize profits
exceeding those that a perfectly competitive market would allow.
Second, to market regulators, who are interested in identifying market
power abuses and developing policies to increase efficiency and social
welfare.

There is an immense amount of literature on strategic bidding in the
context of electricity markets (see e.g., \cite{li2011modeling} for a
related non-exhaustive literature review) particularly when one takes
into account the specific market rules that apply. Currently,
U.S. markets involve multi-part bids for energy and commitment costs, as
well as several types of ancillary services, resulting in
location-dependent hourly and real-time (5 min) prices. In the European
day-ahead market coupling problem \cite{Euphemia}, even more complex
bids are allowed. Regardless of the underlying framework and market
rules, an optimal bidding strategy aims to answer the same question: how
to bid in order to maximize profits.

From a game-theoretic point of view, the approaches for equilibrium
analysis of the strategic bidding problem can be further classified as
Bertrand models, Cournot and Stackelberg models, and {\em Supply
  Function Equilibrium (SFE)} models. In the latter approach, instead of
setting their price bids (Bertrand) or quantities (Cournot), suppliers
bid their supply functions that link prices with quantities. The SFE
literature originates from the seminal work of Klemperer and
Meyer~\cite{klemperer1989supply}, and, since its first application in
electricity markets~\cite{green1992competition}, it has been extensively
studied --- an overview is presented in \cite{holmberg2010supply}, in
both stylized examples and actual electricity markets (see, e.g.,
\cite{anderson2002using, anderson2008finding, holmberg2009numerical} for
analytical and numerical results, and \cite{sioshansi2007good} for an
empirical analysis).

One of the main criticisms regarding game-theoretic approaches is the
unrealistic assumption that the payoff functions of all participants are
publicly available. Most related works deal with this issue by assuming
some type of uncertainty. An early work~\cite{david1993competitive}
proposes a recursive dynamic programming approach for determining the
optimal bid price for each block of generation, in which each supplier
models the uncertainty about rival bid prices by a probability
distribution. In~\cite{wen2001strategic}, the developed bidding scheme
maximizes the hourly profit assuming all other producers' bids are
represented by a multivariate normal distribution whose parameters are
estimated from historical data. In \cite{li2005strategic}, each supplier
assumes types (based on the cost structure) of other suppliers and their
joint probability distribution, based on the published information on
fuel contracts, availability of transmission lines, and operating
parameters. In \cite{yucekaya2009strategic}, a decomposition-based
particle swarm optimization method is proposed to solve the expected
profit maximization problem with the market clearing price modeled as an
uncertain, exogenous variable. In \cite{cherukuri2016decentralized}, a
decentralized Nash equilibrium learning strategy is presented in a
Bertrand competition framework to solve the economic dispatch
problem. Recently, in \cite{mitridati2017bayesian}, a Bayesian inference
approach is proposed to reveal the aggregate supply curve in a day-ahead
electricity market.  In~\cite{motalleb2017non}, a non-cooperative game
with incomplete information among demand response aggregators is
considered under different market conditions, where a Bayesian approach
is used to estimate the unknown information such as the types of
competitors.  In~\cite{zou2016pool}, a multi-period market equilibrium
problem is considered to study the strategic behavior of energy storage
systems, where the optimality conditions of all participants' profit
maximization problems are collected and solved together.

In this paper, we consider SFE-based equilibrium strategies for
suppliers in the context of a day-ahead electricity market.  We address
the aforementioned criticism by estimating payoff functions using an
{\em inverse optimization} approach combined with the theory of
variational inequalities \cite{begupa-vis-mathprogr}.  Inverse
optimization seeks to recover input data to optimization problems from
optimal solutions; it was first introduced in~\cite{ahuja2001inverse}
but recently revisited in new settings~\cite{begupa-BL-12,
  begupa-vis-mathprogr, zhao-pas-seg-gbio-16}.  Interestingly, inverse
optimization has not been extensively used in the context of electricity
markets. In \cite{ruiz2013revealing}, inverse optimization is used to
identify the bids of marginal suppliers in a multi-period
network-constrained electricity pool. In \cite{saez2016data}, it is
employed to address the market-bidding problem of a cluster of
price-responsive consumers of electricity. In \cite{birge2017inverse},
it is used to determine market structure from commodity and
transportation prices; the methodology is applied to data from the MISO
electricity market. Recently, \cite{lu2018data} used inverse
optimization to estimate how loads respond to demand response price
signals.

A preliminary conference version of this paper has appeared
in~\cite{electricity-inverse-equili-cdc2017}. That paper focused on
comparisons between inverse optimization-based strategic bidding versus
earlier approaches~\cite{wen2001strategic}.  The present paper uses a
different, more realistic parametrization of the unknown cost function
with respect to observable market variables, has an extensive numerical
validation of the proposed approach, establishes new rigorous
  results on algorithm termination, and offers extensions to
location-dependent prices.

To the best of our knowledge, we are the first to leverage inverse
optimization for estimating cost parameters in payoff functions and
obtaining equilibrium bids in the context of electricity markets.  Our
method is driven purely by data, in the sense that only the observed
samples are utilized for inference and estimation, without relying on
any distributional assumptions on the observed data. We note that any
hypothesis on the data generating pattern could be questionable, due to
the lack of supporting evidence on such assumption and the noisy nature
of the data. By contrast, data-driven approaches receive input from the
observed samples and are self-adjusted in the estimation process as more
samples are available. 

The main idea of this paper is to learn from past bids of
  electricity suppliers that bid strategically in a day-ahead market
  context. We develop an inverse optimization approach for estimating
  suppliers' cost functions, based on historical bidding data. We
  propose an algorithm that randomly searches for cost function
  parameters with good out-of-sample performance, among multiple
  possible values that are compatible with past data. Our proposed
  framework is validated with extensive numerical experimentation, and
  is extended to accommodate location-dependent prices.

The remainder of the paper is organized as follows. In Section~\ref{s2},
we present the general market framework. In Section~\ref{forward}, we
formulate the strategic bidding problem (referred to as the ``forward''
problem), and in Section~\ref{inverse} we present the inverse
optimization framework as it applies to our day-ahead market setting. In
Section~\ref{s4}, we discuss the specific algorithm we use to estimate
competitors' cost parameters, based on which equilibrium bids are
obtained, and establish its convergence properties. In Section~\ref{s5}
we illustrate our approach with numerical examples. In Section~\ref{s6},
we discuss the extension of our approach to location-dependent
prices. We conclude and provide directions for further research in
Section~\ref{s7}.

\section{Market Framework} \label{s2}

We consider a day-ahead electricity market setting, which is composed of
$N$ electricity suppliers, and a market operator instantiating a {\em
  Power Exchange (PX)}. Each supplier submits a bid curve (supply curve)
that describes the relationship between energy price and production
quantity, for each of the 24 hours of the next day. After receiving the
bidding functions from all suppliers, the market operator clears the
market by balancing aggregate supply and demand; the output is the
hourly market clearing price and the supplier specific dispatch
schedules. Assuming no inter-temporal coupling in the PX setting, the
auctions for different hours are performed separately and
independently. This allows us to consider the bidding strategy for a
specific hour and omit the time index in our analysis.

In actual power markets, the bidding functions are piecewise-constant
curves, reflecting the constant bid price (marginal as-bid cost) for
each block of electricity generation. These piecewise-constant curves
correspond to piecewise-linear functions for the total as-bid costs of
the suppliers, which approximate a quadratic cost function of typical
generators. Piecewise-linear functions are used in practice to allow
reliance on available commercial optimization solvers (for solving
large-scale security-constrained unit commitment and economic dispatch
problems typically formulated as mixed integer linear programming
problems). In this paper, we assume the same affine bid curve as
commonly used in the SFE literature; this assumption not only
facilitates our analysis, but also corresponds to a quadratic
approximation of generator cost functions.

Assume that supplier $i$ submits a linear non-decreasing bid function to
the market operator, $\alpha_i+\beta_i P_i$, $i=1,\ldots, N,$ that
denotes the marginal as-bid cost of power at production level $P_i$, and
$\alpha_i,$ $\beta_i$ are the bidding coefficients to be determined
under the optimal bidding strategy. After receiving these linear bidding
functions, the market operator derives the clearing price and the
generator dispatch schedule as follows:
\begin{align} \label{PX}
		R&=\alpha_i+\beta_iP_i, \ i\in \scrI,\\
		Q&=Q^{\text{forecast}} - \textstyle\sum_{i\in \underline{\scrI}}
                P_i^{\text{min}} - \textstyle\sum_{i\in \overline{\scrI}}
                P_i^{\text{max}}, \notag\\ 
                Q&=\textstyle\sum_{i\in \scrI} P_i, \notag\\
                \overline{\scrI}&=\{
i:\ R \geq \alpha_i + \beta_i P_i^{\text{max}}\},\ \underline{\scrI}=\{
i:\ R \leq \alpha_i + \beta_i P_i^{\text{min}}\}, \notag\\
\scrI&=\{1,\ldots,N\}\setminus \{ \overline{\scrI} \cup \underline{\scrI}\}, 
\notag 
\end{align}
where $R$ is the market clearing price, $Q^{\text{forecast}}$ is the
demand forecast (or as-bid load as is also the case in actual
electricity markets) that is publicly announced by the market operator,
$P_i^{\text{min}}, P_i^{\text{max}}$ are the minimum and maximum
generation levels, respectively, of supplier $i$, $\overline{\scrI}$ is the
set of suppliers producing at $P_i^{\text{max}}$, $\underline{\scrI}$
is the set of suppliers producing at $P_i^{\text{min}}$, $\scrI$ the
set of marginal suppliers, and $Q$ the effective demand met by marginal
suppliers. 
  
Since for $i\in \scrI$ the capacity constraints are not binding, for a
given $Q^{\text{forecast}}$ (hence, $Q$) the solution to (\ref{PX})
becomes
\begin{align}  \label{soltoPX}
		R(\balpha, \bbeta; Q)&=\frac{Q+\sum_{i\in \scrI}
                  \frac{\alpha_i}{\beta_i}}{\sum_{i\in \scrI}
                  \frac{1}{\beta_i}},\\ 
		P_i(\balpha, \bbeta; Q)&=\frac{R(\balpha, \bbeta;
                  Q)-\alpha_i}{\beta_i}, \ i\in \scrI,\notag \\
           P_i(\balpha, \bbeta; Q)&= P_i^{\text{min}},\ i\in
           \underline{\scrI},\; P_i(\balpha, \bbeta; Q)= P_i^{\text{max}},\ i\in
           \overline{\scrI}, \notag 
\end{align} 
where we write $R(\balpha, \bbeta; Q), P_i(\balpha, \bbeta; Q)$ to
explicitly express the dependency on $\balpha \triangleq (\alpha_1,
\ldots, \alpha_N)$ and $\bbeta \triangleq (\beta_1, \ldots, \beta_N)$
for a given demand forecast resulting in effective demand $Q$. 

\section{Forward Problem Formulation} \label{forward} 

The forward problem deals with the individual profit maximization
problem, in which supplier $i$ determines her bidding curve $(\alpha_i,
\beta_i)$ to maximize her profit $\phi_i(\balpha, \bbeta;
Q)$ defined as 
\[ 
\phi_i(\balpha, \bbeta;
Q)=R(\balpha, \bbeta; Q)P_i(\balpha, \bbeta, Q)-C_i(P_i(\balpha, \bbeta;
Q)),
\] 
where $C_i(P_i(\balpha, \bbeta; Q))$ is the production cost at
generation level $P_i(\balpha, \bbeta; Q)$.  The problem is formulated
as follows:
\begin{equation} \label{promax}
\begin{array}{rl} 
 \max\limits_{\alpha_i, \beta_i} & \quad \phi_i(\balpha, \bbeta; Q)\\
 \text{s.t.}   & \quad 0 \le \alpha_i \le \bar{\alpha},\\
& \quad \beta_i > 0,
\end{array}
\end{equation}
where $\bar{\alpha}$ is an upper bound on $\alpha_i$, and is related
  to the price cap in electricity markets.~\footnote{Strictly
    speaking, the price cap imposes a bound on the marginal cost at the
    maximum capacity, $\alpha_i + \beta_i P_i^{\text{max}} \leq
    \bar{\alpha}$, implying that the upper bound is different for each
    supplier $i$, i.e., $\bar{\alpha_i} = \bar{\alpha} - \beta_i
    P_i^{\text{max}}$. However, in practice, the price cap is high
    enough, and we can assume without loss of generality a common upper
    bound on $\alpha_i$.}

Note that the form of the profit function is generic (defined as
  revenues minus costs) and, in general, each supplier can have her own
  cost function. For the case of electricity generators, a common
assumption is a quadratic cost function:
\begin{equation} \label{quadf}
C_i(P_i)= c_{i0}+c_{i1}P_i+ c_{i2} P_i^2, 
\end{equation}
which implies a marginal cost equal to $c_{i1} + 2 c_{i2} P_i$. Since
the intercept of the quadratic cost function (parameter $c_{i0}$) is
constant (this practically refers to the so-called ``no-load'' cost of a
generator), its value will not affect the profit maximization problem of
the supplier and, without loss of generality, it can be set to zero.

A direct comparison of the marginal cost ($c_{i1} + 2 c_{i2} P_i$) and the linear bid function ($\alpha_i + \beta_i P_i$) indicates that the cost parameters $c_{i1}$ and $2 c_{i2}$ correspond to the bidding curve parameters $\alpha_i$ and $\beta_i$, respectively. Hence, truthful bidding
would result in $\alpha_i = c_{i1}$ and $\beta_i = 2 c_{i2}$.

In this paper, we assume that suppliers game only with parameter
$\alpha_i$, and that $\beta_i$ is small and equal to $2 c_{i2}$,
representing a publicly known, technology-specific efficiency decline
associated with approaching generating capacity. 
This assumption corresponds to the ``bid-$\alpha$'' game in
\cite{hu2007using}, implying that $\beta_i$ is known to
other suppliers for all intents and purposes. 
We elaborate on this assumption next.
 
The technology and capacity of individual generating plants is public
information that provides useful partial information about their cost
functions. Nevertheless, their fuel and variable maintenance cost, and
their exact heat rate (efficiency), reflected primarily in parameter
$c_{i1}$, is proprietary and not known with sufficient accuracy to
competitors so as to allow them to bid optimally.  On the other hand, as
also pointed out in \cite{hobbs2000strategic}, the marginal cost
functions of individual suppliers usually have very shallow slopes, and
thus $\beta_i$ is relatively small (about two orders of magnitude lower than $\alpha_i$); furthermore, if both $\alpha_i$ and
$\beta_i$ can be chosen, the existence of a unique equilibrium is rare.
Hence, one can argue that the small value of the slope of the marginal
cost ($2 c_{i2}$) is more or less known, and that the suppliers reflect
this cost in parameter $\beta_i$, \footnote{The parameter $\beta_i$
  reflects $(i)$ the smoothing/regularization of the bid conforming to
  the monotonically increasing market rule requirement (marginal costs
  are physically not strictly monotonic), and $(ii)$ the advantage of
  and desire for achieving unique price-directed marginal generator
  schedules. In the experiments, the values of $\beta_i$ are on the order of
    $0.1$.} as in \cite{hobbs2000strategic}, thus avoiding bidding 
non credible high slopes.  Still, in our results section, we mainly
explore cases in which we allow errors in the estimates of past bids or
the knowledge of parameters (including parameter $\beta_i$).
Furthermore, we note that if we consider the framework from the
perspective of a regulator, the technology-specific data (e.g.,
heat-rate curves) that mainly drive the slope of the marginal cost are
declared by the participants; as such, the slope is relatively easily
calculated.

For the purposes of this paper, we assume that $c_{i1}$ consists of two
cost components: $(i)$ one non-fuel cost component that reflects
operational and maintenance variable costs (e.g., labor, parts, consumables, lubricants, chemicals, consumption from power station supplies, etc.), and $(ii)$ a
fuel-cost component (essentially depending on the heat rate and the
fuel price).  As such, $c_{i1}$ is defined as
\begin{equation} \label{ci1}
 c_{i1} = \theta_{i1} + \theta_{i2} \xi,
\end{equation}
where $\xi$ is a variable reflecting the publicly known fuel price, and
$\theta_{i1}, \theta_{i2}$ are the unknown cost coefficients.
This decomposition is in line with the declared characteristics
of the generation units, which comprise the heat rate curve 
and operational (other than fuel) and maintenance variable costs.
We note that the framework can support even more detailed decompositions
(e.g., consider separately a carbon price for emissions).
Also, we note that the unknown cost coefficients can be interpreted
in various ways. For instance, assuming a publicly known fuel price $\xi$,
coefficient $\theta_{i2}$ may contain the combined effect of the heat rate curve and potential discounts that suppliers may have secured;
such contract information is not available to either regulators or competitors.

Following the above assumptions, i.e., setting $c_{i0} =0$ and $c_{i2} =
(1/2) \beta_i$ in the quadratic cost function (\ref{quadf}), and
using (\ref{ci1}), the profit function for supplier $i$ can be rewritten
as
\begin{align} \label{profit}
\phi_i(\balpha, \bbeta; Q, \xi)
 = & [ R(\balpha, \bbeta; Q)  - (\theta_{i1} + \theta_{i2} \xi ) ]
 P_i(\balpha, \bbeta; Q) \notag \\
&  - \frac{1}{2} \beta_i [ P_i(\balpha, \bbeta; Q) ]^2.
\end{align}
Note that for a given $Q$ and $\xi$, the profit of supplier $i$ is
determined by the actions $\balpha$ of all players and her own cost
parameter $\btheta_i \triangleq ( \theta_{i1}, \theta_{i2})$. 
We therefore write $\phi_i(\btheta_i ; \balpha, Q, \xi)$ to emphasize
this dependency ($\bbeta$ is removed since it is constant and known). 

Since all suppliers choose their bids by solving the profit maximization
problem (\ref{promax}), we can construct a SFE model describing the
game among all profit-maximizing suppliers. By definition, a specific
$\balpha$ is a Nash equilibrium if no single supplier can increase her
profit by unilaterally changing her own bid. We know (see
\cite{electricity-inverse-equili-cdc2017}) that there exists a unique
Nash equilibrium $\balpha^*=(\alpha_1^*, \ldots, \alpha_N^*)$ in this
SFE model, since $\phi_i(\btheta_i ; \balpha, Q, \xi)$ is strictly
concave in $\alpha_i$, i.e., its second partial derivative is strictly
negative.

Next, we compute both the first and the second partial derivatives of
the profit function with respect to $\alpha_i$. From (\ref{profit}), the
first derivative is
\begin{multline} \label{der1}
\nabla_i\phi_i(\btheta_i ; \balpha, Q, \xi) = \\
 \frac{1}{\beta_i} [ \hat \beta_i Q_i 
   + \alpha_i ( {\hat \beta_i }^2 -1 )
 -(\hat \beta_i -1)(\theta_{i1} + \theta_{i2} \xi ) ],
\end{multline}
where $0 < \hat \beta_i = (1/\beta_i)/ \sum_{l\in \scrI} (1/\beta_l) < 1 $, and
$Q_i = ( Q +\sum_{k\in \scrI, k \neq i} \alpha_k/\beta_k)/ \sum_{l\in
  \scrI} (1/\beta_l)$. We 
observe that the first derivative is linear in $\alpha_i$ and also
linear in $\theta_{i1}$ and $\theta_{i2}$. From (\ref{der1}), the second
derivative is
\begin{equation} \label{der2}
\nabla_{ii}^2\phi_i(\btheta_i ; \balpha, Q, \xi) = \frac{1}{\beta_i} (\hat \beta_i - 1 )(\hat \beta_i + 1 ).
\end{equation} 
From (\ref{der2}), it is easy to verify that
$\nabla_{ii}^2\phi_i(\btheta_i ; \balpha, Q, \xi) < 0$, which
implies strict concavity since $\beta_i >0$ and $ 0 < \hat \beta_i
< 1$. 

\section{Inverse Problem Formulation} \label{inverse}

The inverse problem seeks to estimate rivals' cost parameters. This
knowledge is required for estimating the objective of the profit
maximization problem (\ref{promax}). The main theoretical foundation is
attributed to \cite{begupa-vis-mathprogr}, where the authors estimate
the utility functions of the players in a Nash game from the observed
equilibrium.

In the context of this paper, we are given (or we can obtain/estimate)
$M$ past equilibrium bids (observations) $\balpha^j = (\alpha_i^j;\ i\in
\scrI^j)$, $j=1,\ldots, M$, where $\scrI^j$ is the set of marginal
suppliers for observation $j$ and is defined similar to $\scrI$ in
(\ref{soltoPX}). The $\balpha^j$ are realized under different residual
demand levels $Q^j$ and fuel prices $\xi^j$, and we are interested in
estimating $\btheta_i$ of supplier $i=1,\ldots,N$. Without loss of
generality, we assume that for each supplier $i$ there are sufficient
observations $j$ at which supplier $i$ was marginal (i.e., $i\in
\scrI^j$), so that there is enough information to estimate $\btheta_i$.
In case this is not true for some suppliers, we can \emph{a priori}
remove them 
from the set $\{1,\ldots,N\}$ (and appropriately adjust (\ref{PX})). For
such suppliers, we simply do not have enough information to estimate
their cost parameters. 
It should be noted, however, that these suppliers will generally correspond
to base loaded units that do not compete in the market. Given the quadratic cost 
function representation, the resulting linear supply curve is associated with a 
broadly construed notion of marginality that will render non-base-loaded units 
marginal during some hours. As long as each unit is marginal in some of the 
observations -- not all units need to be marginal in all observations -- the 
proposed framework is broadly applicable, and there is no loss of generality from 
the exclusion of base loaded units.

The estimates for the cost parameters are
obtained by applying \cite[Theorem 3]{begupa-vis-mathprogr}, which is
derived through duality, and leads to the following optimization
problem:
\begin{equation} \label{back}
\begin{array}{rl}
\min\limits_{\substack{\by,\bepsilon\\  
\btheta_1, \ldots, \btheta_N}} & \|\bepsilon\|_\infty \\
 \text{s.t.} & 
y_i^j\ge 0,\;  i\in \scrI^j;\  j=1,\ldots,M,\\
&y_i^j \ge \nabla_i\phi_i( \btheta_i; \balpha^j, Q^j, \xi^j),\
\forall i\in \scrI^j, j,\\
&\sum\limits_{i\in \scrI^j} \Bigl(\bar{\alpha}y_i^j-\alpha_i^j\nabla_i\phi_i(\btheta_i;
\balpha^j, Q^j, \xi^j)\Bigr)\le \epsilon_j,
\forall j,\\
&\nabla_i\phi_i( \btheta_i; \balpha^{k_i}, Q^{k_i}, \xi^{k_i}) =
{\phi_i}^{\text{norm}},\ \forall i,\\
\end{array}
\end{equation}
where $\by = (y_i^j)_{j=1,\ldots,M}^{i\in \scrI^j}$ is the decision variable (introduced as a dual variable in \cite[Theorem 2]{begupa-vis-mathprogr}); $\bepsilon=(\epsilon_1,\ldots,\epsilon_M)$,
$\|\bepsilon\|_\infty=\max_j |\epsilon_j|$ is the infinity norm, and the
last constraint is used for normalization purposes and will be discussed
below. We note that the variables in $\nabla_i\phi_i$ are $ \btheta_i$,
and that $\nabla_i\phi_i$ is linear in $\btheta_i$, where $\balpha^j,
Q^j,$ and $ \xi^j$ are parameters of the optimization problem. From
(\ref{der1}), we have
\begin{multline} \label{der1a}
\nabla_i\phi_i(\btheta_i; \balpha^j, Q^j, \xi^j) = 
              \frac{\hat \beta_i}{\beta_i} \frac{ Q^j +\sum_{m\in
                  \scrI^j, m \neq
                  i}  \frac{\alpha_m^j}{\beta_m} }{ \sum_{l\in \scrI^j}
                \frac{1}{\beta_l} } + \\ \frac{\alpha_i^j}{\beta_i}  
( {\hat \beta_i }^2 -1 ) 		 
		+ \theta_{i1} \frac{1 - \hat \beta_i}{\beta_i}
		+ \theta_{i2} \frac{1 - \hat \beta_i}{\beta_i} \xi^j.
\end{multline}

Interestingly, we can reformulate the optimization problem (\ref{back})
as a {\em Linear Programming (LP)} problem, which can be solved very
efficiently. Specifically, instead of the infinite norm objective, we
can introduce constraints that impose an upper bound to each
$|\epsilon_j|$ and then minimize this upper bound.

The last constraint in (\ref{back}) is a normalization constraint, which
is equivalent to \cite[Eqs. (39d) and (39e)]{begupa-vis-mathprogr}. The
right hand side (rhs) of the constraint, ${\phi_i}^{\text{norm}}$, is
some estimate of the partial derivative at a specific point, which is
evaluated at an observation $k_i$ (potentially different for each
$i$). In \cite{begupa-vis-mathprogr}, for illustration purposes, the rhs
estimate is obtained using the actual values of $ \btheta_i$ at a median
observation, considering some lower bounds for the bidding
coefficients. We further elaborate on the implementation of this
constraint in Section~\ref{s4}.

We note that this inverse optimization technique still applies even when
more constraints are imposed in the forward problem setup or when the
bid function is changed, as long as $R$ and $P_i$ have closed-form
expressions w.r.t. the bidding coefficients. 

The quality of the computed equilibrium strategies depends on the
explanatory value of the estimated cost parameters. Indeed, good
estimators should explain future equilibria as well as the equilibria
used to estimate them. The following result, which is a restatement of
\cite[Theorem 6]{begupa-vis-mathprogr}, ensures the quality of the
estimated cost functions under mild conditions. To simplify the
notation, we assume that all suppliers are marginal at all past
observations $j$; otherwise, proper adjustments to the statement of the
theorem can be made.
\begin{thm} \label{theta} Suppose that $\balpha^j, \ j=1,\ldots,M$ are
  i.i.d. realizations of a random variable $\tilde{\balpha}$, and
  $\tilde{\balpha} \in \{\balpha: 0 \le \alpha_i \le \bar{\alpha}, \
  \forall \ i\}$ almost surely. Then, for any $0<\delta<1$, with
  probability at least $1-\eta$ w.r.t. the sampling,
\begin{equation} \label{prob}
		\begin{split}
			\mathbb{P}\big(& \tilde{\balpha} \   \text{is a $z$-approximate equilibrium for the game} \\ 
			& \text{with payoffs defined through}\  \hat{\btheta}_1, \ldots, \hat{\btheta}_N \big) \ge 1-\delta, \\		
		\end{split}
\end{equation}
where $\eta = \sum_{i=0}^{2N} \binom{M}{i} \delta^i (1-\delta)^{M-i}$;
$z$ is the optimal value of problem (\ref{back}); and $\hat{\btheta}_1,
\ldots, \hat{\btheta}_N$ are the optimal solutions to (\ref{back}).
\end{thm}  

Roughly speaking, the $z$-approximate equilibrium describes the
situation where each supplier does not necessarily play her best action
given what others are doing, playing instead a strategy that is no worse
than $z$ relative to the best response. For the definition of
$z$-approximate equilibrium, we refer the interested reader to
\cite[Section 2.2]{begupa-vis-mathprogr}.
 
There are two probability measures in the statement of Theorem
\ref{theta}. One is related to the new data $\tilde{\balpha}$, while the
other is related to the samples $\balpha^1, \ldots, \balpha^M$. The
probability in (\ref{prob}) is taken w.r.t. the new data
$\tilde{\balpha}$. For a fixed set of samples, (\ref{prob}) holds with
probability at least $1-\eta$ w.r.t. the measure of samples. Theorem
\ref{theta} essentially states that given typical samples, the
probability that the estimated cost functions explain well a new future
equilibrium is bounded below. It guarantees the accuracy of the
estimated cost parameters under mild conditions.

\section{Algorithmic Implementation} \label{s4}

In this section, we present the algorithmic implementation for
estimating the rivals' profit functions (or cost parameters
$\btheta_i$), which can then be used to obtain equilibrium bids. We use
historical data from which we can derive the past bids. Suppose we are
aware of the market clearing price and the dispatch schedules of all
suppliers.~\footnote{Such information is publicly available in some
European power markets, or it can be assumed to be discoverable at a
later point in time by market participants.
It is certainly available to regulators even in pool-based markets,
and, to some extent, it can be estimated by entities with market knowledge;
as we will discuss later, errors in the estimates can be viewed as 
noise in the data.} 
Using this information,
the past bids $\balpha^j = (\alpha_i^j;\ i\in \scrI^j)$, $j=1,\ldots,
M$, can be computed via the market-clearing condition in (\ref{PX}),
where $\bbeta$ is constant and known. As before, we assume that for each
supplier $i$ there are sufficient observations $j$ at which supplier $i$
was marginal ($i\in \scrI^j$).

It is worth mentioning that (\ref{back}) might give multiple optimal
solutions. Our goal is to recover the true cost parameters from this
set. Although there might be multiple cost function estimates that can
explain the observed equilibria well, only true costs are expected to
have good out-of-sample performance. The following Algorithm is thus
proposed to identify the true cost functions, based on which equilibrium
bids could be computed via an iterative best response process. We refer
to this algorithm as ``random search,'' since it searches randomly in
the set of optimal solutions until the one that performs well on a
validation dataset is found.

In what follows, we present the main steps of the Algorithm. We define $d$ to be the average discrepancy between computed and true bids on the validation dataset, which evaluates the out-of-sample performance of our solution. The parameter $k$ counts the total number of iterations (random searches). The Algorithm stops when the discrepancy $d$ is small enough or the total number of iterations exceeds a limit.

\hrulefill
\begin{algorithmic}[1]
  \State \textbf{Input:} $N$ suppliers, with constant and known bidding
  slopes $\beta_i$, $i = 1, \ldots, N$; $M$ past bids (observations),
  and for each bid $j= 1,\ldots, M$, the market-clearing price $R^j$,
  residual demand $Q^j$, fuel price $\xi^j$, dispatch schedules $P_i^j$,
  upper bound for bids $\bar{\alpha}$; percentage of training samples
  $p$; tolerance level $\tau$; maximum number of iterations $MaxIter$.

  \State \textbf{Initialize:} $ d = \infty$, $k = 0$. 

  \While{$ d \ge \tau$ and $k < MaxIter$}

  \State $k \leftarrow k + 1$. Randomly choose $M_t=\floor*{Mp}$ samples from all past bids
  (observations) to constitute the training dataset (as a percentage $p$
  of the entire set), and use the remaining bids ($M_v=M-M_t$) as the
  validation dataset.  

  \State Obtain $\hat{\btheta}_i, \ i=1,\ldots,N$, by solving problem
  (\ref{back}) using the training dataset.  

  \State Compute equilibrium strategies (solving (\ref{promax}) via an
  iterative best response process) $\hat{\balpha}_{\text{val}}^j, \
  j=1,\ldots, M_v$, on the validation dataset using $\hat{\btheta}_i, \
  i=1,\ldots,N$.

  \State Evaluate the discrepancy between computed and true bids
  on the validation dataset as
  \begin{equation} \label{c} d = \frac{ \sum_ {j=1} ^ {M_v}
      {\|\balpha_{\text{val}}^j -
        \hat{\balpha}_{\text{val}}^j\|_1}/N}{M_v},
\end{equation}
where $\balpha_{\text{val}}^j$ is the $j$-th true bid (obtained form the
historical data) on the validation dataset, and $\|\cdotp\|_1$ is the
$\ell_1$ norm operator defined as the sum of the absolute elements of
the argument.
\EndWhile

\State Compute equilibrium bids using $\hat{\btheta}_i, \ i=1,\ldots,N$
for given $Q$ and $\xi$.
\end{algorithmic}
\hrulefill

The iterative best response process mentioned in step 6 for computing 
equilibrium bids (which also applies to step 9) is a standard fixed point 
iteration process. Each supplier solves problem (3) assuming all other 
suppliers are fixed in their previous bids (in fact, problem (3) in our 
case can be solved even analytically). Then the bids are updated and the 
process is repeated until an equilibrium is reached, i.e., no supplier can 
gain by unilaterally changing her bid. In practice, this process 
terminates in a few iterations since the profit functions are strictly 
concave.
  
We also note that the algorithm is amenable to parallelization, as
essentially, given adequate resources, all iterations (steps 3 to 8)
could be run in parallel. Our next result establishes that the
  algorithm requires more than $T$ iterations with a probability that
  diminishes exponentially with $T$. Equivalently, we can select a large
  enough maximum number of iterations, $MaxIter$, so that the algorithm
  will terminate before $MaxIter$ is reached with a desirable large
  probability. The result further establishes that the convergence rate
  of the algorithm improves as we increase the size $M_t$ of the
  training set. The proof of the result is included in the Appendix.  We
  numerically explore in Section \ref{s5} the out-of-sample performance
  of the cost estimators obtained through this algorithm.
\begin{thm} \label{thm:rate} Assume that for some $\gamma>0$,
  $(1-\hat{\beta}_i^2)/\beta_i \ge \gamma, \forall i=1, \ldots, N$, and
  the conditions of Thm.~\ref{theta} hold. Assuming that all the past
  bids are at most $\bar{\epsilon}$ away from the equilibrium, and for a
  threshold $\tau = \sqrt{\bar{\epsilon}/(N\gamma)}$, it follows:
\begin{enumerate}
\item for any $T \ge 1$, the probability that the algorithm terminates
  after $T$ iterations is no more than $\eta^T$;
\item for any $0<\epsilon<1$, when $T \ge (\log \epsilon)/(\log \eta),$ the
  probability that the algorithm terminates after $T$ iterations is no
  more than $\epsilon$,
\end{enumerate}
where $\eta$ is defined in the statement of Thm.~\ref{theta}. Moreover,
as we increase the training sample size $M_t$, the number of iterations
that are needed for termination decreases when $M_t$ is large enough.
\end{thm}

Another issue mentioned in the previous section is the implementation of
the normalization constraint, i.e., the last constraint in
(\ref{back}). For the purposes of this paper, unlike
\cite{begupa-vis-mathprogr}, we do not use the true costs (since they
are indeed unknown); instead, we set the rhs (estimate of the partial
derivative) to zero for the median observation of the training dataset.

Finally, we note that the algorithm can handle cases in which ``noise''
is present in the data, e.g., in the past bids (observations). This is
perhaps the most interesting --- and not trivial, application which we
also explore in Section \ref{s5}.

\section{Numerical Illustration} \label{s5}

In this section, we use synthetic input data to test the validity of our
approach. We first describe the experimental setup.

We consider a setup with $N = 2, 3, 4, 5,$ and $10$ suppliers. For each
case, we assume that the true cost parameters $\theta_{i1}$ and
$\theta_{i2}$, as well as $c_{i2}$ of supplier $i$ are equispaced in the
intervals $[7, 5]$, $[0.7, 0.9]$, and $[0.05, 0.07]$,
respectively.~\footnote{For instance, for the case $N = 3$, we have for
  supplier 1, $\theta_{11} = 7$, $\theta_{12} = 0.7$, $c_{12} = 0.05$,
  for supplier 2, $\theta_{21} = 6$, $\theta_{22} = 0.8$, $c_{22} =
  0.06$, and for supplier 3, $\theta_{31} = 5$, $\theta_{32} = 0.9$,
  $c_{32} = 0.07$.}

We generate $M = 200$ past observations, in which demand $Q$ and fuel
price $\xi$ are randomly selected within the intervals $[50, 100]$, and
$[10, 30]$, respectively. For each demand and fuel price realization,
i.e., for each observation $j$ (among the 200), using the true cost
estimators, we generate equilibrium bids $\balpha^j$, assuming that all
suppliers are marginal at all observations. The upper bound
$\bar{\alpha}$ is set to 200.~\footnote{Note that, in practice, we would
  estimate these bids, using the market outcomes.}  The training and
validation datasets are assumed to be of equal size, $M_t = M_v = 100$,
using $p=0.5$. For evaluating the out-of-sample performance we generate
a test dataset with $100$ additional demand and fuel price values,
randomly selected within the aforementioned intervals.

The algorithm was implemented using Matlab R2017a and Gurobi 7.5.1 (for
solving the LP problem (\ref{back})), without any parallelization, and
the computational experiments were run on an Intel i7 5500U, at 2.4GHz,
with 8 GB RAM.

In what follows, we consider two setups: a ``clean'' setup without noise
in the data (Subsection~\ref{s5a}), and a setup with noise
(Subsection~\ref{s5b}). We evaluate the out-of-sample performance
  for the noisy data case (in Subsection~\ref{s5c}), and we perform
  sensitivity analysis with respect to key parameters (in
  Subsection~\ref{s5d}). Lastly, we present a interesting comparison of our approach with a method introduced in \cite{wen2001strategic} (Subsection \ref{wendavid}).

\subsection{``Clean'' Setup} \label{s5a} 

\begin{figure}[hbt]
\centering
\includegraphics[width=2.5in]{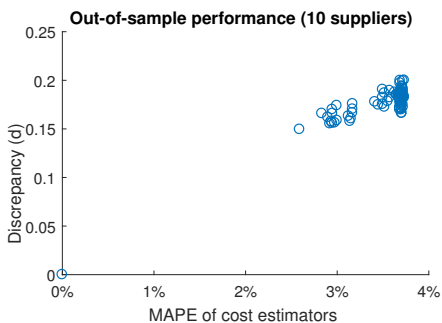}
\vspace{-10pt}
\caption{Validation dataset performance ($N=10$ suppliers, ``clean'' setup).}
\label{fig1}
\end{figure}

As a measure of error for the cost estimators of $\btheta_i$, we use the
{\em Mean Absolute Percentage Error (MAPE)}, defined as
$(100/2N)\sum_{i=1}^N\sum_{l=1}^2
|(\theta_{il}-\hat{\theta}_{il})/\theta_{il}|$. Notably, we expect to
see that cost parameters with low discrepancy values would exhibit low
MAPE values. In all cases, the algorithm managed to exactly reveal the
true costs within 1 or 2 iterations for $N=2,3,4,5$ and 86 iterations
for $N=10$. We plot the results for $N=10$, in Fig. \ref{fig1}. Each
circle corresponds to the discrepancy and MAPE values on the validation
dataset for a certain partition of the samples into training/validation
datasets.

In Fig.~\ref{fig1} we see that in general, the lower the 
discrepancy (the better the performance in the validation set), the
lower the MAPE of the cost estimators (the better the estimates of the
true costs).

\subsection{Setup with Noise} \label{s5b}

We introduce noise in past bid observations, by generating them within
1\% of the optimal value obtained by the iterative best response
process. Alternatively, one may think of this noise as an error in
estimating past bids from the market results. 
For instance, noisy data may be due to errors in estimates of the
market outcome when not all data required (e.g., the exact schedules)
are available, or due to errors in the slopes of the marginal cost
function (reflected in parameter $\beta_i$).
In a more loose interpretation of approximate equilibria, 
one may also think of this noise as observations
in which suppliers do not play exactly their equilibrium bids.

For practical purposes, we set $MaxIter = 10,000$ and a tight tolerance
level $\tau = 10^{-3}$.  In all cases the iteration limit is reached
first; computational times ranged from 5 to 10 min. We note that by
  selecting such a tight tolerance level in absolute figures (with an
  average bid $\alpha$ of around 20, the algorithm would terminate when
  reaching discrepancies lower than 0.005\%), and with 1\% noise present
  in the data, it is almost certain that the algorithm will terminate by
  reaching the iteration limit.  Hence, it is highly unlikely that such
  tight tolerance limits (which are achieved in the clean setup) will be
  also achievable in the case of noisy data. Nevertheless, we keep both
  termination conditions in the framework for the sake of completeness
  in case of noise-free data, to avoid unnecessary iterations. We also
  note that, by setting a relatively high iteration limit, we enhance
  the confidence in our results (see also Thm.~\ref{thm:rate}); we
  elaborate further on the selection of this limit in Subsection
  \ref{s5d}.

\begin{figure}[hbt]
\centering
\includegraphics[width=2.8in]{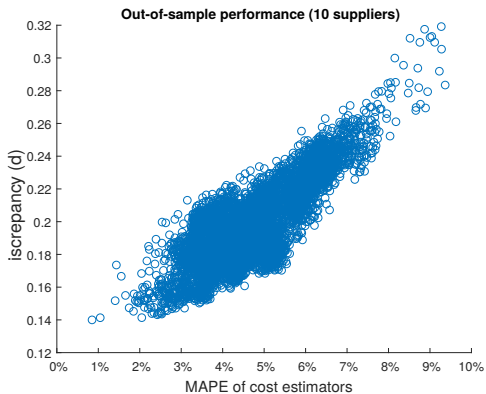}
\vspace{-10pt}
\caption{Validation dataset performance ($N=10$ suppliers, setup with noise).}
\label{fig2}
\end{figure}

For all cases ($N=2,\ldots,5,10$), the best achieved discrepancy at the
validation dataset ranges from 0.111 to 0.154 in absolute figures --- an
average bid $\alpha$ of around 20 implies that the discrepancy is
less than 1\%. The results indicate reasonably good cost estimators with
MAPE ranging from 0.86\% to 3.44\% (for the aforementioned best achieved
discrepancies). In Fig.~\ref{fig2}, we show the performance in the
validation dataset for $N = 10$ suppliers --- compare with
Fig.~\ref{fig1} for the ``clean'' setup.
Fig.~\ref{fig2} also verifies the expected behavior when noise is
present in the data, i.e., good performance in the validation dataset is
associated with good cost estimators.

\begin{figure}[htb]
\centering
\includegraphics[width=3.4in]{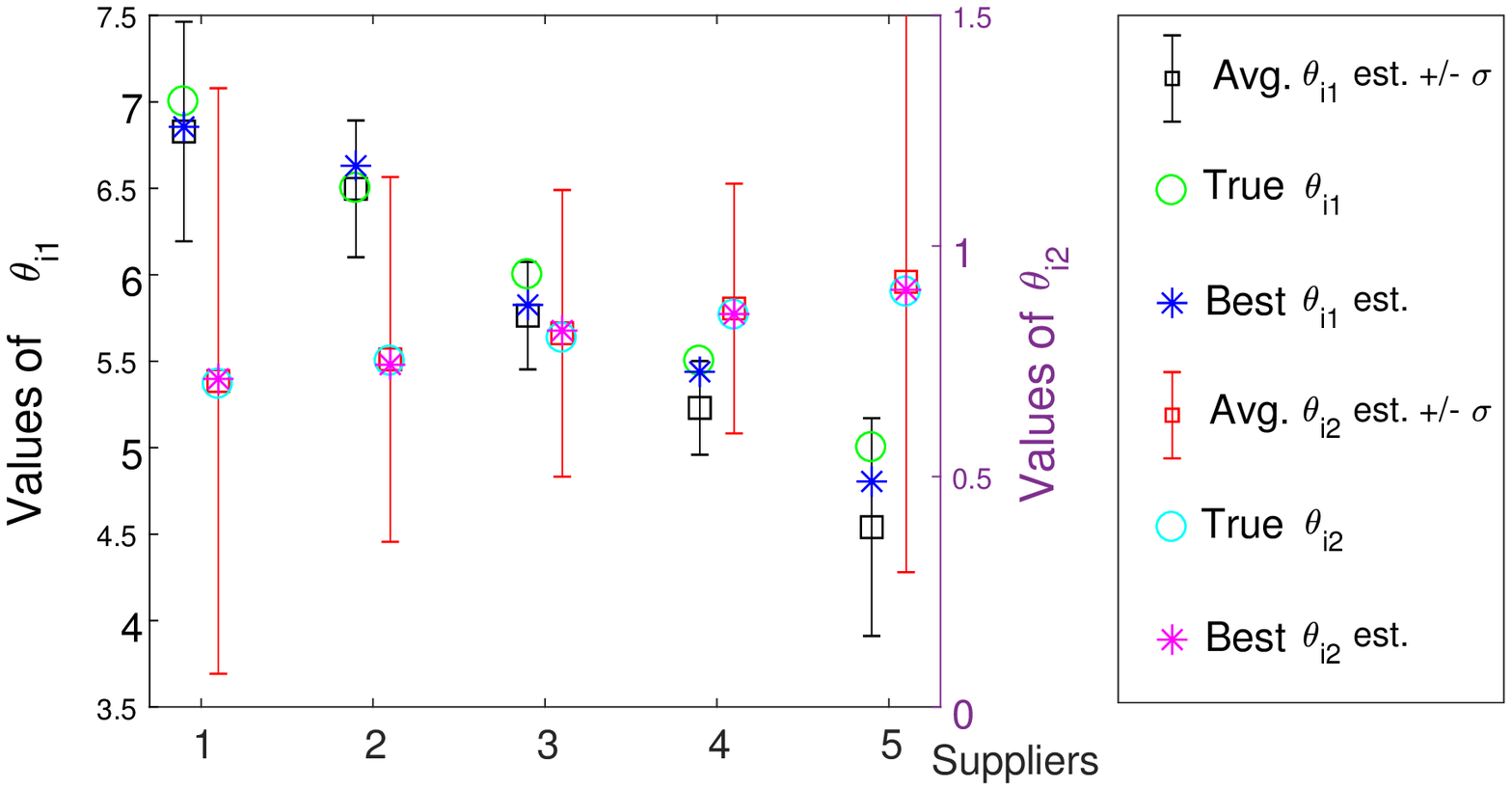}
\includegraphics[width=3.4in]{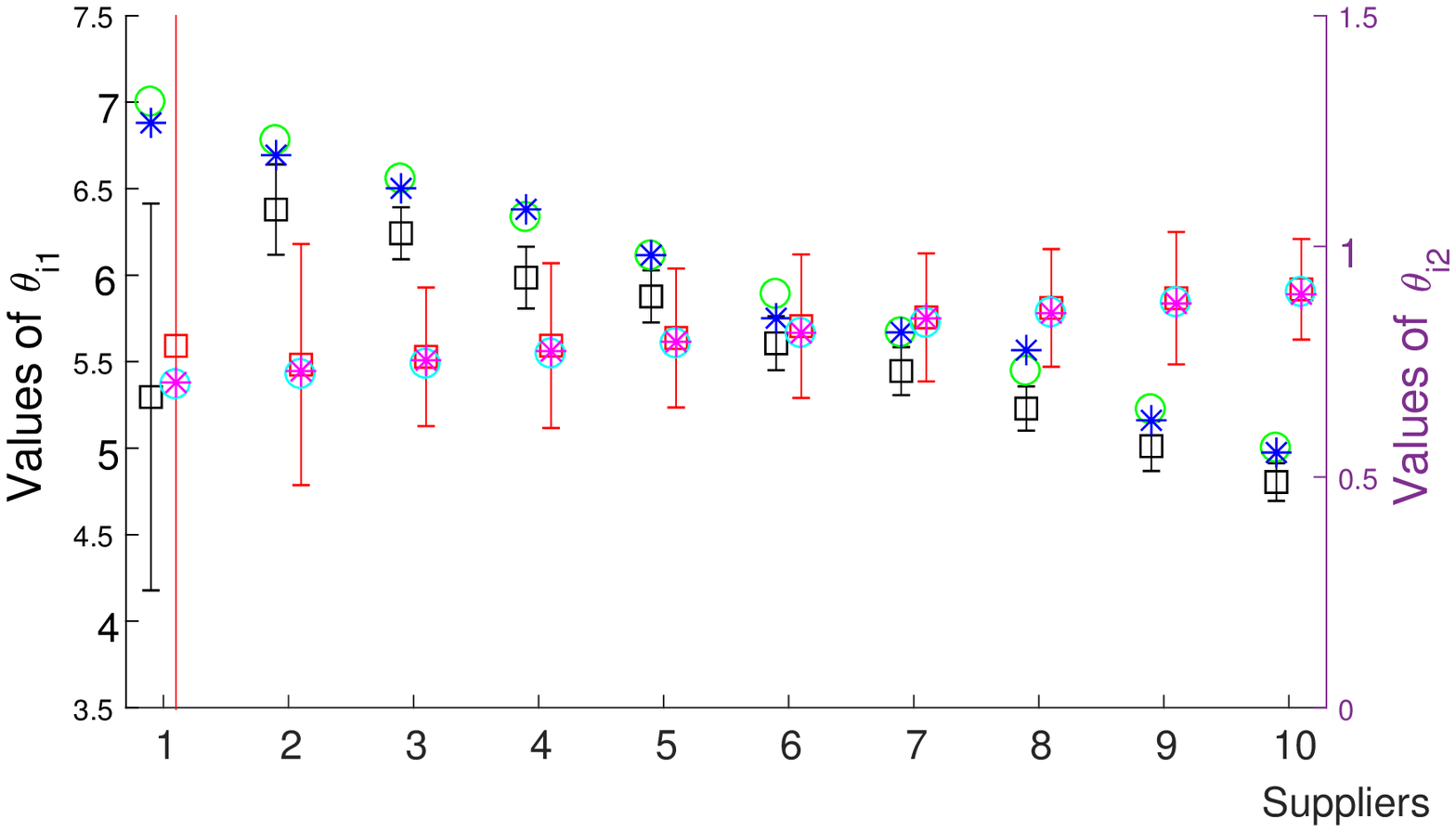}
\vspace{-10pt}
\caption{True costs and estimates (best and average +/-$\sigma$) for $N=5$ (upper figure) and $N=10$ (lower figure) suppliers. Values of $\theta_{i1}$ are shown on the left axis; values of $\theta_{i2}$ are shown on the right axis.}
\label{fig3}
\end{figure}

In Fig.~\ref{fig3}, we plot the true costs $\theta_{i1}$
and $\theta_{i2}$ and their estimates for $N=5$ and
$N=10$ suppliers. We illustrate the best estimate (the
one that corresponds to the lowest discrepancy calculated at the
validation dataset), as well as the average cost estimators and their
standard deviation ($\sigma$) over the 10,000 iterations.

The results indicate that the best estimates range from $-3.9$\%
($-2.3$\%) to $2$\% ($2.2$\%) for $\theta_{i1}$, and from $-1$\%
($-0.4$\%) to $2$\% ($1.3$\%) for $\theta_{i2}$ for the case of $N=5$
($N=10$) suppliers. 
The best estimates, which the Algorithm is designed to obtain, are
reasonably close to the true cost parameters, and hence, they are
expected to exhibit good out-of-sample performance. In fact, even the
average estimates over the 10,000 iterations are not too far from the
true cost parameters.

Next, we evaluate the out-of-sample performance of the best cost
estimators using the test dataset (100 observations, different from the
200 observations that were used for the training and validation
datasets).
\subsection{Out-of-Sample Performance} \label{s5c}
\begin{figure}[hbt]
\centering
\includegraphics[width=3in]{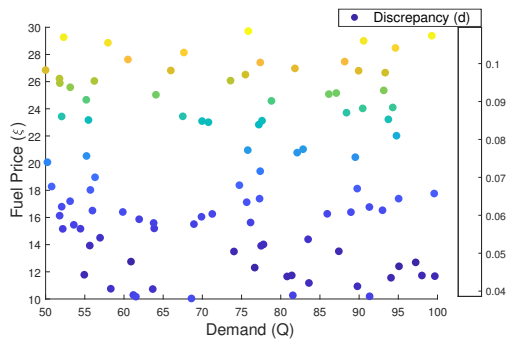}
\vspace{-10pt}
\caption{Discrepancy of bids using the best cost estimators over the
  test dataset ($N=5$ suppliers).}
\label{fig4}
\end{figure}

For each observation (i.e., value of $Q$ and $\xi$) of the
test dataset, we compare the equilibrium bids using the estimates of
$\hat \btheta_i$ with bids derived using the true costs $ \btheta_i$,
and we calculate the discrepancy ($d$) --- using (\ref{c}) with $M_v =
1$. We summarize the results (average discrepancy --- in absolute
figures, and its standard deviation) in Table \ref{tab1}. In
Fig. \ref{fig4}, we illustrate the values of the discrepancy for each
observation of the test dataset for the case of $N=5$ suppliers.
\begin{table}[ht] 
	\caption{Out-of-Sample Performance of Best Cost
          Estimators}  \label{tab1} 
\vspace{-10pt}
	\begin{center}
		\begin{tabular}{c c c }
			\hline
		Suppliers & Avg. Discrepancy ($d$) & Std of $d$ \\ \hline
			2  & 0.086 & 0.0555 \\ 
			3  & 0.047 & 0.0228 \\
			4  & 0.052 & 0.0167  \\
			5  & 0.063 & 0.0220  \\
			10 & 0.104 & 0.0218  \\
			\hline
		\end{tabular}
	\end{center}
\end{table}

The results in Table \ref{tab1} and Fig. \ref{fig4} indicate a
satisfactory out-of-sample performance. In fact, the average discrepancy
is lower than the best achieved over the validation set. The reason is
that the validation dataset contains noise, whereas the test dataset is
a ``clean,'' free of noise setup.

We then take a closer look at the equilibrium bids and profits. We
consider 3 instances of demand and fuel price: (a) $Q = 45$, $\xi = 8$,
(b) $Q = 75$, $\xi = 20$, and (c) $Q = 110$, $\xi = 35$. Note that
instance (b) refers to the mean demand and mean fuel price of the past
bids, whereas instances (a) and (c) contain values that are outside the
intervals used for generating the past bids, i.e., that fall outside the
range of previously observed values. For each instance, we list the
discrepancy values (in absolute figures) in Table \ref{tab2}.
\begin{table}[ht] 
	\caption{Discrepancy of Equilibrium Bids (Estimates vs True
          Costs)}  \label{tab2} 
\vspace{-10pt}
	\begin{center}
		\begin{tabular}{c c c c}
			\hline
          & Instance (a) & Instance (b) & Instance (c) \\		Suppliers & $Q = 45, \xi=8$ & $Q = 75, \xi = 20$ &$Q = 110, \xi = 35$ \\ \hline
			2  & 0.154 & 0.037 & 0.278  \\ 
			3  & 0.048 & 0.041 & 0.114  \\
			4  & 0.104 & 0.035 & 0.111  \\
			5  & 0.062 & 0.059 & 0.140  \\
			10 & 0.064 & 0.100 & 0.162  \\
			\hline
		\end{tabular}
	\end{center}
\end{table}

\begin{figure}[ht]
\centering
\includegraphics[width=3.4in]{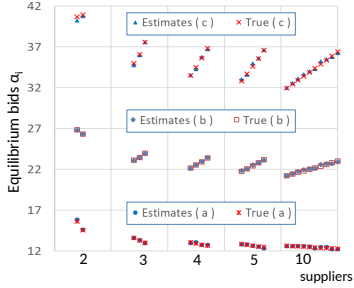}
\vspace{-10pt}
\caption{Equilibrium bids using true costs and estimates for instances (a), (b) and (c) and $N$ = 2, 3, 4, 5, and 10 suppliers (setup with noise).}
\label{fig5new}
\end{figure}

We show the equilibrium bids ($\balpha_i$)
for the three instances, and for $N$ = 2, 3, 4, 5 and 10 suppliers, in Fig. \ref{fig5new}.

From Fig. \ref{fig5new} we see that the difference between the equilibrium
bids using the estimates compared to the ones using the true costs is
very small. For example, for the case $N = 5$ suppliers, the difference ranges from $-0.135$ to 0.047 for instance (a), from
$-0.074$ to 0.133 for instance (b), and from $-0.076$ to 0.351 for
instance (c). The differences reported as percentages range from
$-1.1$\% to 0.4\%, from $-0.3$\% to 0.6\%, and from $-0.2$\% to 1\%, for
instances (a), (b), and (c), respectively.

Lastly, we provide in Table \ref{tab3} the total profits for the three
instances, calculated using (\ref{profit}); the profits using the
estimated costs are shown first and the profits using the true costs
follow in parenthesis for comparison. Apart from the expected behavior
of total profits decreasing when the number of supplier increases, and
total profits increasing with increasing demand, the results also show
that the differences using the estimated costs are very small (in fact the relative differences are within 5\%).
\begin{table}[htb] 
	\caption{Total Profits using Estimates (vs True Costs)}  \label{tab3}
\vspace{-10pt}
	\begin{center}
		\begin{tabular}{c c c c}
			\hline
          & Instance (a) & Instance (b) & Instance (c) \\		Suppliers & $Q = 45, \xi=8$ & $Q = 75, \xi = 20$ &$Q = 110, \xi = 35$ \\ \hline
2	&188.8 (181.7)	&516.0	(518.7)	& 1122.5 (1151.0)\\
3	&81.3 (80.4)	&230.6 (233.7)	& 525.7	(537.5)\\
4	&51.5 (50.3)	&150.8 (150.0)	& 358.4	(361.1)\\
5	&34.5 (36.4)	&113.5	(111.7)	& 294.9	(283.6)\\
10	&15.1 (15.3)	&61.4	(58.3)	& 206.0	(195.8)\\
			\hline
		\end{tabular}
	\end{center}
\end{table}

%
\subsection{Sensitivity Analysis} \label{s5d}
In this subsection, we perform sensitivity analysis with respect to the level of noise, the number of available observations, and the number of iterations. As a base case, we consider the case for $N=5$ suppliers.

\subsubsection{Noise Level} 
\begin{figure}[ht]
\centering
\includegraphics[width=3.4in]{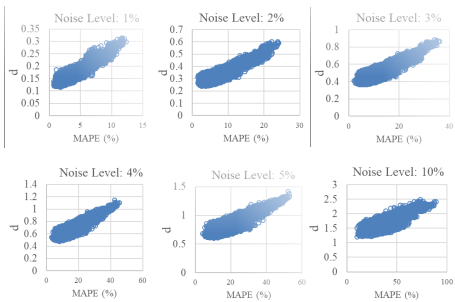}
\vspace{-10pt}
\caption{Validation dataset performance, $N=5$ suppliers, various noise levels.}
\label{fig6}
\end{figure}

As already mentioned, introducing noise in the past bids can be 
thought of as introducing errors in obtaining the past bids from the 
available or estimated market data.
In the previous subsections, we assumed a noise level of 1\%;
in this subsection, we explore higher noise levels (2\%, 3\%, 4\%, 5\%, and 10\%), 
and we present the results for the performance on the validation dataset
(discrepancy vs. MAPE) in Fig. \ref{fig6}.

The best discrepancy achieved ranges from 0.116 (noise level 1\%) 
to 1.204 (noise level 10\%). MAPE values range from 1.75\% to 10.46\%.
The results verify the expected behavior, i.e., discrepancy and
MAPE increase with the noise level. They also verify that  good 
performance in the validation dataset (low discrepancies) is associated 
with good estimators (low MAPE), under all noise levels. Table \ref{tab4} 
presents the out-of-sample performance under various noise levels.
The results are in good agreement with the ones presented in 
Table \ref{tab1}. They also indicate low average discrepancies, which 
increase with the noise level. Note that at high noise levels (see e.g., 10\%)
  the standard deviation becomes low, indicating that the discrepancy is mostly
  affected by the noise. 

\begin{table}[ht] 
	\caption{Out-of-Sample Performance of Best Cost
          Estimators, $N=5$ Suppliers, Various Noise Levels}  \label{tab4} 
\vspace{-10pt}
	\begin{center}
		\begin{tabular}{c c c } 
			\hline
			Noise Level & Avg. Discrepancy ($d$) & Std of $d$ \\ \hline
			1\%  & 0.063 & 0.0220 \\ 
			2\%  & 0.152 & 0.0373 \\
			3\%  & 0.219 & 0.0802  \\
			4\%  & 0.303 & 0.0895  \\
			5\%  & 0.380 & 0.1210  \\	
			10\% & 0.614 & 0.0181  \\ 
			\hline 
		\end{tabular}
	\end{center}
\end{table}

\subsubsection{Past Observations}
\begin{figure}[ht]
\centering
\includegraphics[width=2.3in]{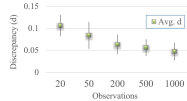}
\vspace{-10pt}
\caption{Out-of-sample performance of best cost estimators (average discrepancy +/-$\sigma$), $N=5$ suppliers, noise level 1\%, various numbers of past observations.}
\label{fig7}
\end{figure}

In the previous subsections, we considered a dataset with $M=200$
available past observations. 
In this subsection, we explore the market outcome for various numbers
of past observations (20, 50, 100, 200, 500, and 1000),
for $N=5$ suppliers, and noise level 1\%. Fig. \ref{fig7} illustrates
the out-of-sample performance of the best cost estimators.

Not surprisingly, and as predicted by Thm.~\ref{theta},
Fig. \ref{fig7} shows that the out-of-sample performance improves with
the number of available past observations. From Thm.~\ref{theta} we see that as the number of past observations $M$ increases, the probability that our estimates $\hat{\btheta}_1, \ldots, \hat{\btheta}_N$ yield an equilibrium that is close to the true one increases, since $\eta$ decreases (this can be seen from the proof in the Appendix). This implies that increasing the training sample size could lead to a small discrepancy between computed and true bids, and thus an improved out-of-sample performance, which is consistent with the observation in Fig. \ref{fig7}. Interestingly, even with 20
available observations, the performance discrepancy is reasonable. We
also checked the average discrepancies (out-of-sample performance) for
the various noise levels with only 20 observations and the results are
also good; average discrepancies range from 0.107 to 0.852 (increasing
with the level of noise) with standard deviations that range from 0.0241
to 0.0644.

\subsubsection{Number of Iterations}
Last but not least, we elaborate on the selection of the maximum number
of iterations, i.e., the termination condition. Thm.~\ref{thm:rate}
yielded a rigorous result that provides guidance on how to select the
iteration limit.  Here, we numerically verify that a small number of
iterations is sufficient to provide some good enough estimators. 
\begin{figure}[ht]
\centering
\includegraphics[width=3in]{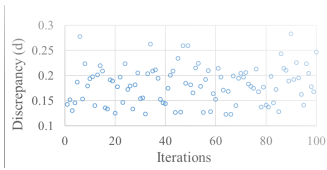}
\vspace{-10pt}
\caption{Discrepancy values calculated in the validation dataset
(step 7 of the algorithm) at each iteration (shown for the first 100 iterations), $N=5$ suppliers, noise level 1\%.}
\label{fig8}
\end{figure}

In Fig. \ref{fig8}, we show the values of the discrepancy 
calculated in the validation dataset at each iteration $k$,
for the case of $N=5$ suppliers with noise level 1\%. For ease of exposition, we plotted the first 100 iterations of the algorithm (out of the 10,000). The results indicate that good enough estimators (with low discrepancy) can be obtained early in the process. In fact, in the first 100 iterations, we observed 7 instances with discrepancies that are less than 10\% higher than the best achieved (which in our tests was 0.116). We elaborate on this indication in the following figure.  
\begin{figure}[ht]
\centering
\includegraphics[width=3in]{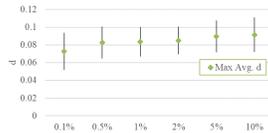}
\vspace{-10pt}
\caption{Out-of-sample performance of cost estimators (maximum average discrepancy +/-$\sigma$), $N=5$ suppliers, noise level 1\%, top x\% (x-axis) of iterations (ranked in ascending order of discrepancy).}
\label{fig9}
\end{figure}

Fig. \ref{fig9} shows the maximum average discrepancy (+/- $\sigma$)
computed over the top x\% of the iterations ranked in ascending order of
discrepancy computed over the same test dataset. For instance, for the
top 1\%, i.e., the top 100 iterations out of 10,000, the worst achieved
out-of-sample performance was an average $d=0.100$ with $\sigma =
0.0167$. 

\begin{figure}[ht]
	\centering
	\includegraphics[width=3in]{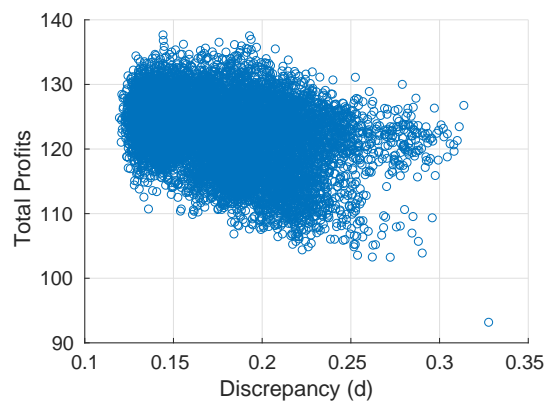}
	\vspace{-10pt}
	\caption{Total profits vs. discrepancy values calculated in the validation dataset, $N=5$ suppliers, noise level 1\%, 200 observations, 10,000 iterations.}
	\label{fig10}
\end{figure}

Lastly, we investigate the relationship between discrepancy values and the total profits. In Fig. \ref{fig10}, we show the values of the discrepancy and the total profits calculated in the validation dataset at each iteration, for the case of $N=5$ suppliers, $M=200$ past observations with noise level $1\%$, and $MaxIter=10,000$ iterations. It can be seen that the larger the discrepancy is, the smaller the total profits, which validates the use of the discrepancy as a performance metric in identifying the best cost estimators.

\subsection{Comparison with the Method in \cite{wen2001strategic}} \label{wendavid}
In this subsection, we compare the performance of our approach with the method introduced in \cite{wen2001strategic}, which assumes
that the bidding coefficients of any rival supplier follow a bivariate normal distribution whose mean and covariance could be inferred from the past observed bids. The profit maximization problem for a single supplier $i$ is then reformulated to involve only $\alpha_i$ and $\beta_i$, with the market clearing price $R$ and the production level $P_i$ evaluated at the mean bids of the rival suppliers.  

We note that \cite{wen2001strategic} does not estimate rivals' cost functions; it calculates the bids using as input the mean and covariance of rivals' bidding coefficients. In our setup, since we assume that the bidding slope $\beta$ is known, only the mean and variance of $\alpha$ need to be inferred from the historical data. For the needs of our comparisons, we use the noisy setup with 1\% noise level.


We evaluate the performance of \cite{wen2001strategic} on the test dataset (that contains 100 equilibrium bids). The average discrepancy between the solution obtained using \cite{wen2001strategic} and the equilibrium bids on the test dataset is shown in Table \ref{tab5} (the differences with the values presented in Table \ref{tab1} that are obtained using our approach are shown in parentheses). The results show that \cite{wen2001strategic} leads to significantly higher discrepancies (as well as higher standard deviations), indicating that our approach obtains a more accurate prediction of the bids.

\begin{table}[ht] 
	\caption{Average Discrepancy of Bids on the Test Dataset Using \cite{wen2001strategic} (Differences with Table \ref{tab1})}  \label{tab5} 
	\vspace{-10pt}
	\begin{center}
		\begin{tabular}{c c c }
			\hline
			Suppliers & Avg. Discrepancy ($d$)  & Std. of  $d$ \\ \hline
			2  & 1.406 (+1.32)  & 0.7891 (+0.73) \\ 
			3  & 1.026 (+0.98)  & 0.5709 (+0.55) \\
			4  & 0.818 (+0.77)  & 0.4551 (+0.44)  \\
			5  & 0.681 (+0.62)  & 0.3793 (+0.36)\\
			10 & 0.372 (+0.27) & 0.2071 (+0.19) \\
			\hline
		\end{tabular}
	\end{center}
\end{table}

%

We then take a closer look at the bids and profits by considering the 3 instances
illustrated in Section \ref{s5c}. In Table \ref{tab7} we list the discrepancy of bids
obtained using \cite{wen2001strategic} (in parentheses we show the differences with
the values listed in Table \ref{tab2}). The results show that for instance (b), which
represents instances that have been constantly observed in the past,
\cite{wen2001strategic} achieves a comparable performance with our approach, and
occasionally even lower discrepancies. However, the performance of
\cite{wen2001strategic} is much worse for instances (a) and (c), which contain values
that were outside the range of past observations. The reason is that
\cite{wen2001strategic} implicitly assumes that rivals' bidding behavior is similar
to what has been observed in the past, which results in a large bias when a new,
unseen scenario occurs. By contrast, through estimating rivals' cost functions from
the past bids, our approach acquires more information, and exhibits a stronger
out-of-sample inference capability that guarantees a low estimation bias for every
possible scenario.  Interestingly but not surprisingly, \cite{wen2001strategic}
achieves better results when the number of suppliers is large, in which case the bids
are very close (due to competition), see e.g., Fig.~\ref{fig5new}.
\begin{table}[ht] 
	\caption{Discrepancy of Bids Using \cite{wen2001strategic} (Differences with Table \ref{tab2})}  \label{tab7} 
	\vspace{-10pt}
	\begin{center}
		\begin{tabular}{c c c c}
			\hline
			& Instance (a) & Instance (b) & Instance (c) \\		Suppliers & $Q = 45, \xi=8$ & $Q = 75, \xi = 20$ &$Q = 110, \xi = 35$ \\ \hline
			2  & 3.773 (+3.62) & 0.048 (+ 0.01) & 4.776 (+4.50) \\ 
			3  & 2.487 (+2.44) & 0.050 (+0.01) & 3.211 (+3.10) \\
			4  & 1.923 (+1.82) & 0.043 (+0.01) & 2.497 (+2.39) \\
			5  & 1.581 (+1.52) & 0.037 (-0.02) & 2.057 (+1.92)  \\
			10 & 0.849 (+0.79) & 0.021 (-0.08) & 1.108 (+0.95) \\
			\hline
		\end{tabular}
	\end{center}
\end{table}

In Table \ref{tab8} we list the total profits for the bids derived using
\cite{wen2001strategic}, and we show in parentheses the differences with the values
in Table \ref{tab3} obtained using the cost estimates. Recall that the profits
obtained using the cost estimates are very close to the profits obtained using the
true costs (within +/- 5\%). For instance (b), the two methods obtain similar
profits. For instance (c) which represents higher demands and fuel prices, our
approach achieves much higher profits. The opposite is true for instance (a), which
represents lower demands and fuel prices. This is an interesting result, which is
explained by the fact that the bids generated using \cite{wen2001strategic} are
biased by the mean past bids. Hence, they tend to be values of past observations,
i.e., they are inflated for the instances of lower demands and fuel prices, and they
are reduced for higher demands and fuel prices. More specifically, the values of the
equilibrium bids are about 25\% higher in instance (a), and about 12\% lower in
instance (c) compared to the ones shown in Fig. \ref{fig5new}. Of course, in instance
(a) where the bids (and hence, profits) are inflated, the conditions are ripe for a
new supplier to come in, underbid, and capture significant market share.

\begin{table}[htb] 
	\caption{Total Profits of Estimated Bids Using \cite{wen2001strategic} (Differences with Table \ref{tab3} Using Estimates)}  \label{tab8}
	\vspace{-10pt}
	\begin{center}
		\begin{tabular}{c c c c}
			\hline
			& Instance (a) & Instance (b) & Instance (c) \\		Suppliers & $Q = 45, \xi=8$ & $Q = 75, \xi = 20$ &$Q = 110, \xi = 35$ \\ \hline
			2	& 350.2 (+161.40)	& 515.3 (-0.70)	& 636.0 (-486.50)\\
			3	& 192.4 (+111.10)	& 230.0 (-0.60)	& 186.2 (-339.50) \\
			4	& 137.1 (+85.60)	& 146.7 (-4.10)	& 87.0 (-271.40)\\
			5	& 107.7 (+73.20)	& 108.9 (-4.60)  & 57.4 (-237.50)\\
			10	& 53.7 (+38.60)	& 56.8 (-4.60)	& 73.7 (-132.30)\\
			\hline
		\end{tabular}
	\end{center}
\end{table}

In conclusion, our approach possesses a stronger out-of-sample inference capability attributed to the estimation of the cost functions. The method proposed in \cite{wen2001strategic} ignores the interaction among suppliers, assumes a normal distribution and uses only the mean values of the past bids to infer rivals' behavior, which accounts for its unsatisfactory performance in new, unseen scenarios.


\section{Extension to Location-Dependent Prices} \label{s6}  

So far, we assumed competition among suppliers in uncongested
networks. Indeed, several day-ahead electricity markets clear ignoring
network congestion. In the instance of such day-ahead market rules, the
system operator adjusts the generation dispatch to observe line flow
capacity constraints and ensure secure and reliable operation.

In U.S. markets, however, the transmission system representation has
been part of the standard market design for many years, with resulting
``Locational Marginal Prices'' (LMPs) representing the marginal cost at
each node of the transmission system. In practice, without entering in a
detailed analysis of how LMPs are formed, we note that ``price islands''
may characterize clearing prices, differing only slightly to reflect
varying loss factors.

Our method appears to assume away the fact that network connected
markets result in location dependent clearing price differentials driven
by $(i)$ small effects of location-specific line loss contributions, but
also, $(ii)$ significant contributions during network congestion
events. It can capture and address, however, significant
congestion-caused differentials by detecting market-splitting
occurrences that result in ``price islands'' with essentially
homogeneous prices within each island. Although this limits the number
of relevant observations when price islanding occurs, it utilizes the
unusually high or low price events associated with congestion.

Our approach applies to price islands where congestion is a result of
generators outside of the relevant island, who do not set the price and,
hence, are not part of our analysis. In this context, the market
clearing price in each island $s$ is $R_s = \alpha_i + \beta_i P_i$, $i
\in \scrI_s$, where $\scrI_s$ is the set of marginal suppliers in island
$s$. The residual demand in each island is $Q_s = \sum_{i \in \scrI_s}
P_i.$ Similarly to (\ref{soltoPX}), the solution for island $s$, in
terms of price and quantity, is
\begin{equation*}  \label{Sol_isl}
\begin{split}
	R_s(\balpha_s, \bbeta_s; Q_s)&= \frac{Q_s +\sum_{i \in \scrI_s}
          \frac{\alpha_i}{\beta_i} }{\sum_{i \in I_s} \frac{1}{\beta_i}
        },\\ 
	P_i(\balpha_s, \bbeta_s; Q_s)&=\frac{R_s(\balpha_s,
          \bbeta_s; Q_s)-\alpha_i}{\beta_i}, \ i \in \scrI_s, 
\end{split}
\end{equation*}  
where $\balpha_s, \bbeta_s$ refer to suppliers in island $s$. 

In the forward problem, the profit of the suppliers in each island is
therefore straightforwardly defined. Considering the inverse problem, we
note that the past bids may contain both congested and uncongested
instances.  The set of marginal suppliers may be different in each
instance, and furthermore the islands in the congested instances may be
different. But in either case, our approach can handle these different
sets, since for each observation $j$, we can have different price
islands $s \in \scrS_j$, where $\scrS_j$ is the set of price islands for
the $j$-th observation, and for each price island of that observation we
have a set of marginal suppliers denoted by $\scrI_s^j$. Hence, the inverse
optimization problem in (\ref{back}) is applied for each price island
$s$ for observations $j$ within this island, and for the respective set
of marginal suppliers for the specific island. Essentially, the
uncongested case represents one single island (and can still be
described by the above notation).

\section{Conclusions and Further Research} \label{s7}

In this paper, we proposed an inverse optimization method to estimate
electricity suppliers' cost functions in the day-ahead electricity
market based on historical bidding data. The problem of computing
optimal bidding strategies can be seen as an equilibrium computation
problem given the estimated payoff functions. We applied a ``random
search'' algorithm to estimate the cost function parameters of
electricity generators; specifically, the parameters that are
proportional to their generation output. The algorithm essentially seeks
cost function parameters (among multiple possible values compatible with
the past data) which have good out-of-sample generalization
performance. We established strong, exponential-type probabilistic
  convergence guarantees for this algorithm. Extensive numerical
experimentation verifies that one can recover accurate estimates of the
cost function parameters, which, in turn, allows generators to bid with
knowledge of how competitors would respond. Even though we considered a
simple setting involving no congestion or transmission network effects,
we discussed an extension of the methodology to location-dependent
prices.

Regarding future research directions, it would be of interest to develop
non-parametric approaches that do not require to assume a specific
parametric form for the cost functions. Finally, in addition to
estimating competing generators' cost functions, our methodology is
particularly useful in estimating the underlying cost functions of
market participants who bid synthetic or virtual generators
corresponding to contracts with either physical generation owners or a
portfolio of demand-response-capable consumers.


\appendix

\section*{Proof of Thm.~\ref{thm:rate}}

\begin{proof}
  Assume that the optimal solutions to the inverse problem at iteration
  $k$ are $\hat{\btheta}_1^k, \ldots, \hat{\btheta}_N^k$, and the
  optimal value is $z_k$. We will first show that the function
  $f^k(\balpha) \triangleq
  \bigl(-\nabla_1\phi_1(\hat{\btheta}_1^k; \balpha, Q, \xi),
  \ldots, -\nabla_N\phi_N(\hat{\btheta}_N^k; \balpha, Q,
  \xi)\bigr)$ is strongly monotone. For simplicity we suppress the
  dependence of $f$ on $\btheta_i, Q$ and $\xi$. By definition, a
  function $f^k(\balpha)$ is {\em strongly monotone} if $\exists
  \gamma>0$ such that
\begin{equation*}
\bigl(f^k(\balpha_1)-f^k(\balpha_2)\bigr)'(\balpha_1 - \balpha_2) \ge \gamma \|\balpha_1 - \balpha_2\|_2^2, \forall \balpha_1, \balpha_2.
\end{equation*}
Plugging in the formula for $\nabla_i\phi_i(\btheta_i; \balpha, Q,
\xi)$, we have:
\begin{multline*}
  f^k(\balpha_1)-f^k(\balpha_2) =\\
  \Bigl(\frac{1-\hat{\beta}_1^2}{\beta_1} (\alpha_{1,1} - \alpha_{2,
    1}), \ldots, \frac{1-\hat{\beta}_N^2}{\beta_N} (\alpha_{1,N} -
  \alpha_{2, N})\Bigr),
\end{multline*}
where $\alpha_{1,i}, \alpha_{2,i}$ are the $i$-th elements of
$\balpha_1$ and $\balpha_2$, respectively. Using
$(1-\hat{\beta}_i^2)/\beta_i \ge \gamma$, $\forall i$, it follows
\begin{equation*}
\bigl(f^k(\balpha_1)-f^k(\balpha_2)\bigr)'(\balpha_1 - \balpha_2) 
\ge \gamma \|\balpha_1 - \balpha_2\|_2^2.
\end{equation*}
With the strongly monotone function $f^k(\balpha)$, we can use
\cite[Thm. 8]{begupa-vis-mathprogr}, which shows that for any
$0<\delta<1$, with probability at least $1-\eta$ with respect to the
sampling,
\begin{equation*}
  \|\balpha_{\text{val}}^j - \hat{\balpha}_{\text{val}}^j\|_2 \le
  \sqrt{z_k/\gamma},\quad \forall j = 1, \ldots, M_v,
\end{equation*}
where $z_k$ is the optimal value of the inverse optimization problem
(\ref{back}) at iteration $k$. Using the norm inequality
\begin{equation*}
\|\balpha_{\text{val}}^j - \hat{\balpha}_{\text{val}}^j\|_1 \le \sqrt{N} \|\balpha_{\text{val}}^j - \hat{\balpha}_{\text{val}}^j\|_2,
\end{equation*}
we obtain that at iteration $k$, the discrepancy satisfies:
\begin{equation*}
  d \le \sqrt{z_k/(N\gamma)} \le \sqrt{\bar{\epsilon}/(N\gamma)},
\end{equation*}
which yields that
\begin{equation} \label{dbound} \mbb{P} (d \le
  \sqrt{\bar{\epsilon}/(N\gamma)}) \ge 1-\eta.
  \end{equation}
  Since the iterations are independent from each other, setting $p
  \triangleq \mbb{P} (d \le \sqrt{\bar{\epsilon}/(N\gamma)})$, and
  using (\ref{dbound}), we have:
  \begin{equation*}
    \mbb{P}\Bigl(\text{Algorithm terminates after $T$ iterations}\Bigr)
    = (1-p)^T \le \eta^{T}.
  \end{equation*}
  Therefore, for any small $0<\epsilon<1$, if the probability that the
  algorithm terminates after $T$ iterations is below $\epsilon$, we need
  $T \ge (\log \epsilon)/(\log \eta)$.  

  We next show that as the training sample size $M_t$ increases, the
  number of iterations that are needed decreases for a large enough
  $M_t$. First note that
\begin{equation*}
\begin{aligned}
	\eta & = \sum_{i=0}^{2N} \binom{M_t}{i} \delta^i (1-\delta)^{M_t-i} \\
		& = \sum_{i=0}^{2N}\frac{(M_t-i+1)\ldots M_t}{i!}\delta^i (1-\delta)^{M_t-i} \\
		& \le \sum_{i=0}^{2N}\frac{(M_t)^i}{i!}\delta^i (1-\delta)^{M_t-i}.
\end{aligned}
\end{equation*}
Define $h_i(M_t) \triangleq \frac{(M_t)^i}{i!}\delta^i
(1-\delta)^{M_t-i}$, and take its derivative:
\begin{equation*}
  \nabla h_i(M_t) = \frac{(M_t)^{i-1}}{(i-1)!}\delta^i (1-\delta)^{M_t-i} \Bigl(1 + \frac{M_t}{i} \log(1-\delta)\Bigr).
\end{equation*}
We see that for a large enough $M_t$, $\nabla h_i(M_t) <0, \forall i$,
since $\log(1-\delta)<0$. Therefore, as $M_t$ increases, $\eta$
decreases, and the number of iterations, i.e., $(\log \epsilon)/(\log
\eta)$, decreases as well.
\end{proof}


\begin{thebibliography}{10}
	\bibitem{mcc1988spot}
	F.~C. Schweppe, M.~C. Caramanis, R.~D. Tabors, and R.~E. Bohn, \emph{Spot
		pricing of electricity}.\hskip 1em plus 0.5em minus 0.4em\relax Boston, MA:
	Kluwer Academic Publishers, 1988.
	
	\bibitem{david2000strategic}
	A.~K. David and F.~Wen, ``Strategic bidding in competitive electricity markets:
	a literature survey,'' in \emph{Power Engineering Society Summer Meeting,
		2000. IEEE}, vol.~4.\hskip 1em plus 0.5em minus 0.4em\relax IEEE, 2000, pp.
	2168--2173.
	
	\bibitem{li2011modeling}
	G.~Li, J.~Shi, and X.~Qu, ``Modeling methods for {GenCo} bidding strategy
	optimization in the liberalized electricity spot market--a state-of-the-art
	review,'' \emph{Energy}, vol.~36, no.~8, pp. 4686--4700, 2011.
	
	\bibitem{Euphemia}
	(2016) Euphemia public description. {PCR} market coupling algorithm. [Online].
	Available: \url{http://www.nordpoolspot.com}
	
	\bibitem{klemperer1989supply}
	P.~D. Klemperer and M.~A. Meyer, ``Supply function equilibria in oligopoly
	under uncertainty,'' \emph{Econometrica}, vol.~57, pp. 1243--1277, 1989.
	
	\bibitem{green1992competition}
	R.~J. Green and D.~M. Newbery, ``Competition in the {British} electricity spot
	market,'' \emph{Journal of Political Economy}, vol. 100, no.~5, pp. 929--953,
	1992.
	
	\bibitem{holmberg2010supply}
	P.~Holmberg and D.~M. Newbery, ``The supply function equilibrium and its policy
	implications for wholesale electricity auctions,'' \emph{Utilities Policy},
	vol.~18, pp. 209--226, 2010.
	
	\bibitem{anderson2002using}
	E.~J. Anderson and A.~B. Philpott, ``Using supply functions for offering
	generation into an electricity market,'' \emph{Operations Research}, vol.~50,
	no.~3, pp. 477--489, 2002.
	
	\bibitem{anderson2008finding}
	E.~J. Anderson and X.~Hu, ``Finding supply function equilibria with asymmetric
	firms,'' \emph{Operations Research}, vol.~56, no.~53, pp. 697--711, 2008.
	
	\bibitem{holmberg2009numerical}
	P.~Holmberg, ``Numerical calculation of an asymmetric supply function
	equilibrium with capacity constraints,'' \emph{European Journal of
		Operational Research}, vol. 199, no.~1, pp. 285--295, 2009.
	
	\bibitem{sioshansi2007good}
	R.~Sioshansi and S.~Oren, ``How good are supply function equilibrium models: an
	empirical analysis of the {ERCOT} balancing market,'' \emph{Journal of
		Regulatory Economics}, vol.~31, pp. 1--35, 2007.
	
	\bibitem{david1993competitive}
	A.~David, ``Competitive bidding in electricity supply,'' in \emph{IEE
		Proceedings C (Generation, Transmission and Distribution)}, vol. 140,
	no.~5.\hskip 1em plus 0.5em minus 0.4em\relax IET, 1993, pp. 421--426.
	
	\bibitem{wen2001strategic}
	F.~Wen and A.~David, ``Strategic bidding for electricity supply in a day-ahead
	energy market,'' \emph{Electric Power Systems Research}, vol.~59, no.~3, pp.
	197--206, 2001.
	
	\bibitem{li2005strategic}
	T.~Li and M.~Shahidehpour, ``Strategic bidding of transmission-constrained
	{GENCOs} with incomplete information,'' \emph{IEEE Transactions on Power
		Systems}, vol.~20, no.~1, pp. 437--447, 2005.
	
	\bibitem{yucekaya2009strategic}
	A.~D. Yucekaya, J.~Valenzuela, and G.~Dozier, ``Strategic bidding in
	electricity markets using particle swarm optimization,'' \emph{Electric Power
		Systems Research}, vol.~79, no.~2, pp. 335--345, 2009.
	
	\bibitem{cherukuri2016decentralized}
	A.~Cherukuri and J.~Cort{\'e}s, ``Decentralized {Nash} equilibrium learning by
	strategic generators for economic dispatch,'' in \emph{American Control
		Conference (ACC), 2016}.\hskip 1em plus 0.5em minus 0.4em\relax IEEE, 2016,
	pp. 1082--1087.
	
	\bibitem{mitridati2017bayesian}
	L.~Mitridati and P.~Pinson, ``A {B}ayesian inference approach to unveil supply
	curves in electricity markets,'' \emph{IEEE Transactions on Power Systems},
	2017.
	
	\bibitem{motalleb2017non}
	M.~Motalleb and G.~R, ``Non-cooperative game-theoretic model of demand response
	aggregator competition for selling stored energy in storage devices,''
	\emph{Applied Energy}, vol. 202, pp. 581--–596, 2017.
	
	\bibitem{zou2016pool}
	P.~Zou, Q.~Chen, Q.~Xia, G.~He, C.~Kang, and A.~J. Conejo, ``Pool equilibria
	including strategic storage,'' \emph{Applied Energy}, vol. 177, pp.
	260–--270, 2016.
	
	\bibitem{begupa-vis-mathprogr}
	D.~Bertsimas, V.~Gupta, and I.~C. Paschalidis, ``Data-driven estimation in
	equilibrium using inverse optimization,'' \emph{Mathematical Programming},
	vol. 153, no.~2, pp. 595--633, 2015.
	
	\bibitem{ahuja2001inverse}
	R.~Ahuja and J.~Orlin, ``Inverse optimization,'' \emph{Operations Research},
	vol.~49, no.~5, pp. 771--783, 2001.
	
	\bibitem{begupa-BL-12}
	D.~Bertsimas, V.~Gupta, and I.~C. Paschalidis, ``Inverse optimization: A new
	perspective on the black-litterman model,'' \emph{Operations Research},
	vol.~60, no.~6, pp. 1389--1403, 2013.
	
	\bibitem{zhao-pas-seg-gbio-16}
	Q.~Zhao, A.~Stettner, E.~Reznik, I.~C. Paschalidis, and D.~Segr\'{e}, ``Mapping
	the landscape of metabolic goals of a cell,'' \emph{Genome Biology}, vol.~17,
	no.~1, p. 109, 2016.
	
	\bibitem{ruiz2013revealing}
	C.~Ruiz, A.~J. Conejo, and D.~J. Bertsimas, ``Revealing rival marginal offer
	prices via inverse optimization,'' \emph{IEEE Transactions on Power Systems},
	vol.~28, no.~3, pp. 3056--3064, 2013.
	
	\bibitem{saez2016data}
	J.~Saez-Gallego, J.~M. Morales, M.~Zugno, and H.~Madsen, ``A data-driven
	bidding model for a cluster of price-responsive consumers of electricity,''
	\emph{IEEE Transactions on Power Systems}, vol.~31, no.~6, pp. 5001--5011,
	2016.
	
	\bibitem{birge2017inverse}
	J.~R. Birge, A.~Horta\c{c}su, and J.~M. Pavlin, ``Inverse optimization for the
	recovery of market structure from market outcomes: An application to the
	{MISO} electricity market,'' \emph{Operations Research}, vol.~65, no.~34, pp.
	837--855, 2017.
	
	\bibitem{lu2018data}
	T.~Lu, Z.~Wang, J.~Wang, Q.~Ai, and C.~Wang, ``A data-driven {S}tackelberg
	market strategy for demand response-enabled distribution sytems,'' \emph{IEEE
		Transactions on Smart Grid}, 2018.
	
	\bibitem{electricity-inverse-equili-cdc2017}
	R.~Chen, M.~Caramanis, and I.~C. Paschalidis, ``Strategic equilibrium bidding
	for electricity suppliers in a day-ahead market using inverse optimization,''
	in \emph{Proceedings of the 56th IEEE Conference on Decision and Control},
	Melbourne, Australia, December 12--16 2017.
	
	\bibitem{hu2007using}
	X.~Hu and D.~Ralph, ``Using {EPECs} to model bilevel games in restructured
	electricity markets with locational prices,'' \emph{Operations Research},
	vol.~55, no.~5, pp. 809--827, 2007.
	
	\bibitem{hobbs2000strategic}
	B.~F. Hobbs, C.~B. Metzler, and J.-S. Pang, ``Strategic gaming analysis for
	electric power systems: An {MPEC} approach,'' \emph{IEEE Transactions on
		Power Systems}, vol.~15, no.~2, pp. 638--645, 2000.
	
\end{thebibliography}
\end{document}